\documentclass{amsproc}

\usepackage{fullpage,graphicx, amssymb, amsmath, mathrsfs}

\newtheorem{theorem}{Theorem}[section]

\newtheorem{proposition}[theorem]{Proposition}
\newtheorem{conjecture}[theorem]{Conjecture}
\newtheorem{openproblem}[theorem]{Open Problem}
\theoremstyle{definition}
\newtheorem{definition}[theorem]{Definition}
\newtheorem{example}[theorem]{Example}

\theoremstyle{remark}

\numberwithin{equation}{section}

\newcommand{\fagnano}{\mathscr{F}}

\begin{document}

\title{Towards the Koch Snowflake Fractal Billiard: \\ Computer Experiments and Mathematical Conjectures}

\author{Michel L. Lapidus\footnote{The research of the first author was supported in part by the U.S. National Science Foundation Grant DMS-0707524.}}
\address{Department of Mathematics, University of California, Riverside, CA 92521--0135, USA}

\email{lapidus@math.ucr.edu}

\author{Robert G. Niemeyer}
\address{Department of Mathematics, University of California, Riverside, CA 92521--0135, USA}

\email{niemeyer@math.ucr.edu}

\date{June 12, 2009.}


\begin{abstract}
In this paper, we attempt to define and understand the orbits of the Koch snowflake fractal billiard $KS$.   This is a priori a very difficult problem because $\partial(KS)$, the snowflake curve boundary of $KS$, is nowhere differentiable, making it impossible to apply the usual law of reflection at any point of the boundary of the billiard table.  Consequently, we view the prefractal billiards $KS_n$ (naturally approximating $KS$ from the inside) as rational polygonal billiards and examine the corresponding flat surfaces of $KS_n$, denoted by $\mathcal{S}_{KS_n}$.  In order to develop a clearer picture of what may possibly be happening on the billiard $KS$, we simulate billiard trajectories on $KS_n$ (at first, for a fixed $n\geq 0$).  Such computer experiments provide us with a wealth of questions and lead us to formulate conjectures about the existence and the geometric properties of periodic orbits of $KS$ and detail a possible plan on how to prove such conjectures.  


\end{abstract}

\maketitle

\section{Introduction}
\label{sec:introduction}

The long-term goal of this work is to justify the existence and investigate, via computer-aided experiments and suitable mathematical arguments, a prototypical fractal billiard, namely, the Koch snowflake billiard.

\begin{figure}
\begin{center}
\includegraphics{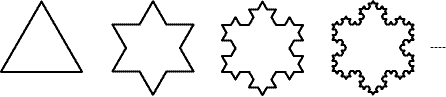} 
\end{center}
\caption{Construction of the Koch snowflake from the equilateral triangle $\Delta$ (with sides of length 1); from left to right: $KS_0=\Delta$, $KS_1$, $KS_2$, $KS_3$. The process continues ad infinitum, yielding the Koch snowflake curve, viewed here as the boundary of the Koch snowflake billiard table $KS$.  For each integer $n\geq 0$, $KS_n$ denotes the $n$th prefractal polygonal approximation to $KS$. Note that $KS_n$ defines a rational billiard because its interior angles are either $\pi/3$ or $4\pi/3$.  While this property is not required to be a rational billiard, the interior angles of $KS_n$ alternate between $\pi/3$ and $4\pi/3$.  (Depending on the context, $KS_n$ also sometimes refers to the polygonal boundary of this rational billiard; the latter boundary is more correctly denoted by $\partial(KS_n)$, however.)}
\label{fig:kochconstruction}
\end{figure}

Since the Koch snowflake curve (the boundary of the Koch snowflake billiard table, see Fig. \ref{fig:kochconstruction}) is nowhere differentiable, it is certainly not clear from the outset that the Koch snowflake billiard---or, let alone, the associated billiard flow---is a legitimate mathematical object of study.  Indeed, the first question that springs to mind is ``\textit{How do you reflect a billiard ball off of a nondifferentiable boundary}?''  On the other hand, because the Koch curve is highly symmetric (indeed, it is \textit{self-similar}) and can be thought of as some kind of `infinite polygon' (see Fig.~\ref{fig:kochconstruction}), it may not be entirely unreasonable to harbor such an expectation.

In this paper, we provide experimental and geometric evidence towards this claim, as well as formulate conjectures and open problems concerning this challenging problem.  In a nutshell, the Koch snowflake billiard is viewed as a suitable limit of (rational) polygonal billiard approximations; and analogously for the associated billiard (and geodesic) flows.  We also identify a variety of (potential) periodic and quasiperiodic orbits of the Koch snowflake fractal billiard, and discuss some of their properties.  

In order to present our experimental results and conjectures about the snowflake billiard, we need to briefly recall several concepts from the theory of dynamical systems (more specifically, of rational polygonal billiards) and fractal geometry.

A \textit{polygonal billiard} $B$ is given by a polygon $P$ and the bounded region enclosed by $P$.  For example, the square billiard is given by the square and the region bounded by the square.  When studying a polygonal billiard, one is primarily concerned with the orbits of a pointmass moving at unit speed in $B$. We assume that our billiard ball experiences no friction and perfectly elastic collisions with the smooth portions of the boundary.  Upon colliding with a smooth portion of the boundary, we reflect at the boundary according to Snell's Law, which says that the angle of incidence equals the angle of reflection; see, e.g., \cite{GaStVo}, \cite{Gu}, Chapter 9 of \cite{KaHa}, or \cite{Ta}. In general, if the billiard ball collides with a corner of the polygonal billiard table $B$, we must terminate the trajectory because we can not determine reflection at such points.\footnote{We will see that there are polygons for which we can determine reflection at a corner, the equilateral triangle billiard being one.}  We call the map that describes the flow on the billiard $P$ the \textit{billiard map}. An important case of a polygonal billiard is a \textit{rational billiard} and is defined to be a polygon $P$ for which every interior angle is a rational multiple of $\pi$; see, e.g., \cite{MaTa}, \cite{Sm}, \cite{Vo} or \cite{Zo}.

Although a fractal is a shape that lacks many of the appealing qualities of polygons, so long as we are dealing with the right fractal, it is a shape that can be generated by polygons, a fact that is heavily exploited in this paper.  The type of fractals we are interested in are self-similar fractal curves in $\mathbb{R}^2$ that are nowhere differentiable and have infinite length.  Specifically, the fractal billiard that we eventually want to define and investigate is the Koch snowflake billiard, denoted by $KS$; the boundary of the associated billiard table is the Koch snowflake curve, denoted by $\partial(KS)$.  The construction of the Koch snowflake fractal curve is given in Fig. \ref{fig:kochconstruction}.  The \textit{prefractal} $KS_n$ of $KS$ is then an $n$th level polygonal approximation of $KS$.  Like many fractals, the Koch snowflake curve $\partial(KS)$ has a non-integer dimension.  In the case of the snowflake, that dimension is $\log_3{4}$, which indicates that the curve $\partial(KS)$ has infinite length.

The Koch snowflake is not a rational polygon.  The boundary of the Koch snowflake `fractal billiard' is nondifferentiable, hence making it very difficult to properly define, and let alone analyze, the billiard map.  However, the behavior of a pointmass may be anything but random, because, as was alluded to just above, $KS$ has the very special property that there is a sequence of finite polygonal approximations $\{KS_n\}_{n=0}^\infty$ converging to $KS$ such that for each finite $n$, $KS_n$ is a rational billiard; see Fig. \ref{fig:kochconstruction}.  In order to overcome the limitations of the Koch snowflake, we may examine the prefractal billiard approximations $KS_n$ and attempt to make an argument in support of the existence of periodic orbits of the limiting fractal billiard $KS$, based on results concerning the rational polygonal billiards $KS_n$.  Therefore, in some sense, the snowflake curve is viewed as an (infinite) `fractal rational polygon' and the associated billiard table $KS$ as a `fractal rational billiard'.

Roughly speaking, a compact set $F\subseteq \mathbb{R}^2$ is said to be \textit{self-similar} if it can be written as a finite union of scaled (or rather, similar) copies of itself.  More specifically, this means that $F$ is a nonempty compact subset of $\mathbb{R}^2$ and there exists a finite collection $\{W_i\}_{i=1}^l$ of contractive similarity transformations of $\mathbb{R}^2$, with $l\geq 2$, such that $F=\bigcup_{i=1}^l W_i(F)$.  It then follows that for any nonempty compact subset $X$ of $\mathbb{R}^2$, we have that $\lim_{k\to\infty} \mathcal{W}^k(X) = F$ (in the sense of the Hausdorff metric), where $\mathcal{W}(X):=\bigcup_{i=1}^l W_i(X)$ and $\mathcal{W}^k$ denotes the $k$th iterate of $\mathcal{W}$; see, e.g., Chapter 9 of \cite{Fa}.  For this reason, the self-similar set $F$ is called the attractor of the iterated function system $\{W_i\}_{i=1}^l$.

\begin{figure}
\begin{center}
\includegraphics{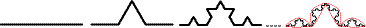}
\caption{Construction of the von Koch curve (left), often simply called the ``Koch curve'' here.  At each stage of the construction, the middle third of each line segment is replaced by the other two sides of an equilateral triangle based on that segment.  The self-similarity of the von Koch curve (right).  The Koch curve is decomposed as the union of four pieces similar to the whole curve.}
\label{fig:kochcurveconstruction}
\end{center}
\end{figure}

\begin{figure}
\begin{center}
\includegraphics{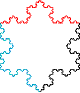}
\caption{The Koch snowflake curve is obtained by pasting together three copies of the von Koch curve (from Fig.~\ref{fig:kochcurveconstruction}).}
\label{fig:3kochcurves}
\end{center}
\end{figure}

The Koch curve is a perfect example of a self-similar set; see Fig.~\ref{fig:kochcurveconstruction}.\footnote{See also [\textbf{Fa}, pp. xiv--xv] and, in a related context, \textbf{[La1,2]}, \cite{LaNRG}, \cite{LaPa} or [\textbf{La-vF}, \S{12.3.1} \& \S{12.5}].}  Strictly speaking, the Koch snowflake fractal is not self-similar; specifically, it is the result of pasting together three copies of the Koch curve; see Fig. ~\ref{fig:3kochcurves}.  In light of this, we make a slight abuse of language by referring to $KS$ (or rather, its boundary $\partial(KS)$) as a self-similar fractal.

We can exploit the fact that the proposed billiard $KS$ has a self-similar boundary.  We notice that, in a sense, $KS_n$ is constructed from equilateral triangles. The billiard ball dynamics on the equilateral triangle billiard $\Delta=KS_0$ are well understood.  In particular, the periodic orbits of the equilateral triangle billiard $\Delta$ are classified in \cite{BaUm}.  The resulting theorem, combined with the fact that $KS_n$ can be embedded in a tiling $T_{\Delta_n}$ of the plane by equilateral triangles with side lengths $1/3^n$, will aid us in explaining our experimental results and formulating some of our conjectures.

This paper is structured as follows. In Section 2, we discuss the necessary theory of rational billiards so that we may better present our results on the prefractal polygonal billiards $KS_n$; see Fig. \ref{fig:kochconstruction}. In \S\ref{subsec:unfolding}, we discuss a useful tool for discerning periodic orbits of rational billiards, namely, the unfolding of billiard trajectories.  In general, such an unfolding gives rise to a surface that depends on the orbit, but we also state known results for constructing a surface $\mathcal{S}_P$ which is independent of the periodic orbit of the rational billiard $P$.  We recall in \S\ref{subsec:surface} that the flow on the rational billiard $P$ corresponds to the geodesic flow on the corresponding surface $\mathcal{S}_P$.  In fact, as is explained in \S\ref{subsec:singularities}, vertices constitute singularities of the billiard map and certain copies of the vertices in $\mathcal{S}_P$ constitute singularities of the flow on the surface.  We  discuss what implications this has for the billiard ball dynamics.

In \S\ref{subsec:KSNratbilliard}, we view the prefractal polygonal billiard $KS_n$ as a rational billiard, and discuss the consequences of this fact for the billiard flow at the vertices of $KS_n$. We also construct the corresponding surface for $KS_1$, $KS_2$ and $KS_3$ and indicate the singularities of each surface.  We show that, in agreement with the general theory of rational billiards, reflection at acute angles can be defined, but that reflection at obtuse angles cannot be determined unambiguously.

In the last part of Section 2, namely \S\ref{subsec:baum}, we discuss the results of \cite{BaUm} on the equilateral triangle billiard and their implications for the billiard $KS_n$.  In doing so, we explain how $KS_n$ can be embedded in a tiling $T_{\Delta_n}$ of the plane by equilateral triangles with side lengths $1/3^n$ (as was alluded to above) and give an equivalence relation on the collection of all periodic orbits of the rational billiard $\Delta$.  We do not seek to generalize the results of \cite{BaUm} to $KS_n$ but instead use these results to provide us with a list of appropriately adjusted initial conditions for testing in $KS_n$.

In Section 3, we present and discuss our experimental results on the prefractal billiard $KS_n$ by raising and sometimes answering a variety of research questions.  In general, we examine the behavior of what we call \textit{induced} orbits of $KS_n$.  An induced orbit is one for which the initial condition was appropriately adapted from an initial condition in the equilateral triangle billiard $\Delta$.  

In \S\ref{subsec:gammanorbits}, we examine the behavior of periodic orbits of $KS_n$ in the collection $\gamma_n$ of billiard paths determined by the initial conditions $(x_{mid},\pi/3)$, where $x_{mid}$ is a midpoint of a side of $KS_n$ and $\pi/3$ is the initial angle of the periodic orbit at the starting point $x_{mid}$.  We attempt to organize the orbits in $\gamma_n$ according to their geometrical or dynamical behavior, and provide definitions for what we call \textit{Fagnano orbits}, \textit{piecewise Fagnano orbits} and \textit{primary piecewise Fagnano orbits} of $KS_n$.  While these are very interesting orbits, the elements in $\gamma_n$ do not constitute all of the induced orbits of $KS_n$.  

In \S\ref{subsec:nongammanorbits}, we examine the behavior of periodic orbits induced by the initial condition $(x_0,\theta_0)$, where $x_0\neq x_{mid}$ or $\theta_0 \neq \pi/3$ and the orbit remains nonsingular, meaning that it does not hit any vertices of the billiard table.  A particularly interesting example of such a periodic orbit is given by $(x_{mid}+\delta x,\pi/6)$, where $\delta x$ is a suitable value (or vector) used to perturb the initial basepoint lying at the midpoint $x_{mid}$ of a side of $KS_n$.  Because of the nature of the equilateral triangle billiard, the reflection at the vertices of $\Delta = KS_0$ can be determined.  Consequently, in \S\ref{subsec:singularorbits}, these singular (and periodic) orbits are then used to induce singular orbits of $KS_n$, for some $n\geq 0$.  

Finally, in \S\ref{subsec:quasiperiodic}, we close Section 3 by discussing a seemingly uninteresting type of orbit, which we call a \textit{quasiperiodic orbit} of $KS_n$.  Such an orbit can be considered as a `rational approximation' to a given periodic orbit.  For example, if $(x_0,\theta_0)$ is an initial condition of a periodic orbit, then a corresponding quasiperiodic orbit would have an initial condition $(x_1,a/b)$, where $a/b$ is a rational approximation of $\theta_0$ (obtained via a continued fraction expansion of $\theta_0$).  Because the boundary of the billiard $KS_n$ is changing as $n$ increases, the study of quasiperiodic orbits of $KS_n$ should enable us to define a suitable notion of `quasiperiodic orbit' of $KS$.

In Section 4, we conclude our paper by stating several conjectures based on our experimental results and provide a list of open problems  and conjectures related to the proposed fractal billiard $KS$ and a particular collection of its periodic orbits.  In particular, we conjecture the existence of what we call the \textit{primary piecewise Fagnano orbit} of $KS$, which we propose would be a suitable limit of primary piecewise Fagnano orbits of the prefractal approximations $KS_n$. In the special case of $pp\fagnano$, the `primary piecewise Fagnano' orbit of $KS$,\footnote{The presumed periodic orbit $pp\fagnano$ of $KS$ is induced by the Fagnano orbit of the original triangle $\Delta = KS_0$ (i.e., the shortest periodic orbit of $\Delta$), along with its appropriate counterpart in each billiard table $KS_n$; see \S\ref{subsec:gammanorbits}, along with Conjectures~\ref{conj:weakppf} and \ref{conj:existenceOfSelfSimilarPeriodicOrbits} of Section 4.} we also conjecture that its `footprint'\footnote{I.e., the subset of the boundary $\partial(KS)$ consisting of all incidence points of the periodic orbit $pp\fagnano$ of $KS$.} on the boundary $\partial(KS)$ is a self-similar subset of $\partial(KS) \subseteq \mathbb{R}^2$; more specifically, it is the natural counterpart of the middle-third Cantor set strung around the Koch snowflake curve $\partial(KS)$.

The other conjectures and open problems stated in Section 4 are clearly of a longer term nature.  They concern, for example, the existence of the Koch snowflake billiard $KS$ (as a proper mathematical object) and of the associated billiard flow,\footnote{as well as of the associated `fractal surface' $\mathcal{S}_{KS}$ and of the corresponding geodesic flow (conjecturally equivalent to the billiard flow on $KS$); see \S\ref{subsec:surface} and \S\ref{subsec:KSNratbilliard}, along with parts (i) and (iii) of Conjecture \ref{conj:fractalbilliardFractalSurface}.} or (in the very long-term) the relationship between the length spectrum of the elusive Koch snowflake fractal billiard and the frequency spectrum of the corresponding Koch snowflake drum (e.g., [\textbf{La1--2},\textbf{LaNRG},\textbf{LaPa}], along with \S{12.3} and \S{12.5} of \cite{La-vF}).\footnote{See Open Problem~\ref{openproblem:fractalbilliardvsfractaldrum}.}

In closing this introduction, we mention that the subject of billiards---and particularly, polygonal billiards or even, rational billiards---is an active area of research in the field of dynamical systems.  Books, survey articles and research papers on various aspects of this topic include, respectively, [\textbf{KaHa},\textbf{Ta}], [\textbf{GaStVo},\textbf{Gu},\textbf{HuSc},\textbf{MaTa},\textbf{Sm},\textbf{Vo},\textbf{Zo}] and [\textbf{BaUm},\textbf{GuJu1--2},\textbf{KaZe},\textbf{Ma},\textbf{Ve1--3}].  To our knowledge, none of these references deals with the elusive case of `fractal billiards'.  It is certainly true, however, that these works (and the many relevant references therein) have played a key role in motivating and guiding our investigations on the `Koch snowflake billiard' (and its rational polygonal approximations) reported in the present paper.

\section{Flows on Rational Billiards and Corresponding Surfaces}

In this section, we want to discuss the properties of the billiard flow on $KS_n$ and describe how $KS_n$ can be viewed naturally as a \textit{rational billiard}.  A \textit{rational polygon} is a polygon with interior angles that are rational multiples of $\pi$. If the boundary of a polygonal billiard is a rational polygon, then we call the associated billiard a \textit{rational billiard}. For us to discuss the flow on the associated closed surface, we need to define a few concepts.\footnote{For an introduction to the theory of rational polygonal billiards from various points of view, we refer, e.g., to \cite{GaStVo}, \cite{Gu}, \cite{MaTa}, \cite{Sm}, \cite{Ta}, \cite{Vo} or \cite{Zo}.}  

\begin{definition}
If $H$ is a regular $N$-gon, where $N$ is a positive integer, then $D_N$ is the group of symmetries of $H$ and has cardinality $2N$.  It is called the \textit{dihedral group} and is generated by the reflections in the lines through the origin that meet at angles $\pi/N$.
\end{definition}

\begin{definition}
Let $A(P)$ be the group of planar motions generated by the reflections in the sides of a polygon $P$.  Furthermore, let $G(P)$ denote the subgroup of the orthogonal group $O(2)$ consisting of the linear parts of the elements of $A(P)$.
\end{definition}

When the polygon $P$ has a connected boundary, then we can give an alternate characterization of it being rational.

\begin{proposition}
Let $P$ be a polygon with connected boundary. Then $P$ is a rational billiard if and only if the group $G(P)$ is finite.  In that case, if the interior angles of the rational polygon $P$ are written in the form $\pi(m_j/n_j)$ for $j=1,...,r$, where $m_j$, $n_j$ are coprime positive integers and $r$ is the number of vertices of $P$, then $G(P)$ is isomorphic to the dihedral group $D_N$, where $N$ is the least common multiple of $n_1,n_2,...,n_r$.
\label{prop:altDefinitionOfRatPolygon}
\end{proposition}

If $P$ is a rational billiard, then we have at our disposal a method for discerning periodic orbits.  This method was first given in  \cite{KaZe} and reiterated in various forms in the literature; see, e.g., the survey articles \cite{GaStVo}, \cite{Gu}, \cite{MaTa} and \cite{Zo}. Let $x_0$ be the initial position of the billiard ball and $\theta_0$ be the initial direction of the billiard ball.  Let $x_1$ be the point on the side $s_1$ of the billiard table at which the billiard ball collides, and $\theta_1$ be the angle at which the billiard ball reflects off of $s_1$. If the orbit $\alpha$ is periodic, then we want to illustrate this by `unfolding' the billiard table and its contents.

\subsection{Unfolding the billiard table $\mathbf{B}$.}
\label{subsec:unfolding}

\begin{figure}
\begin{center}
\includegraphics{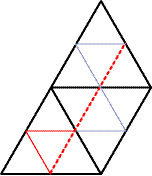}
\caption{Unfolding the equilateral triangle billiard:  the orbit unfolds to a straight line segment (denoted by the dashed line).  The terminal point on the upper right equilateral triangle in the unfolding corresponds to the initial point and the angles are the same.  Hence, this corresponds to a periodic orbit of the equilateral triangle.}
\label{fig:unfoldtheorbit}
\end{center}
\end{figure}

To unfold the billiard table and its contents, reflect the contents of $B$ in the side $s_1$.  This results in a segment continuing in the direction of $\theta_0$.  Moreover, we see that we have created a situation where the trajectory passes through $s_1$ and continues on to a segment that is collinear with the initial segment of the trajectory (see Fig. \ref{fig:unfoldtheorbit} for an example in the equilateral triangle billiard).  Continuing this process, we see that we can unfold a billiard table and its contents.  When the orbit is a collection of line segments, such an unfolding will be a straight line.  Furthermore, the length of the unfolding is exactly the length of the orbit in $B$.  When the orbit unfolds and terminates on an edge of a polygon $\tilde{P}$ in the unfolding such that the terminal point is a copy of the initial point and the segment makes an angle with the side of $\tilde{P}$ equal to the initial angle, then the original orbit of $P$ is periodic.

It is significant that we can unfold an orbit to determine the periodicity of the said orbit.  In order to construct a closed surface out of the unfolding (which is determined by the periodic orbit), we identify sides of the unfolding based on where the ball collides.   In particular, we identify the terminal side with the initial side since the orientation of $\tilde{P}$ and the original polygon  $P$ are the same.  This results in a closed surface\footnote{By `closed', we mean that the surface does not have any boundary.} whose geodesic flow\footnote{which is actually a straight-line flow on the closed surface} corresponds to the billiard flow on $P$. We can see this clearly in Fig.~\ref{fig:unfoldtheorbit}. However, this closed surface depends on the initial position and angle.  We want to construct a surface $\mathcal{S} = \mathcal{S}_P$ that is independent of these initial conditions. In addition, the geodesic flow on that surface $\mathcal{S}_P$ will turn out to be equivalent to the original billiard flow on $P$.

\subsection{The invariant surface $\mathcal{S}_P$ and its associated flow.}
\label{subsec:surface}

Consider the product $P\times A(P)$. Essentially, $P\times A(P)$ contains $2N$ copies of our rational billiard table $P$. A graphical representation of $P\times A(P)$ can be given by the following.  If we fixed a vertex of a polygon $P$ such that reflection in the adjacent sides of $P$ generates $2N$ copies of $P$, then the resulting shape is what is called a \textit{generalized polygon}.  If two sides of the generalized polygon are translates of each other and/or the result of a rotation by $\pi$, then we consider these two sides to be equivalent.  Then, upon modding out by this equivalence relation, we obtain a closed surface $\mathcal{S}_P := P\times A(P)/\sim$ that does not depend on any particular orbit.

\begin{figure}
\begin{center}
\includegraphics{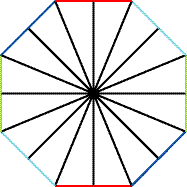}
\caption{The closed surface $\mathcal{S}_P$ corresponding to the triangle $P(3\pi/8,\pi/8,\pi/2)$; the associated generalized polygon is an octagon.  Pairs of appropriate sides should be identified as indicated, in order to obtain $\mathcal{S}_P$, a surface of genus $g=2$.  Note, however, that strictly speaking, this planar representation is not valid.  We illustrate such an embedding here as a heuristic device.}
\label{fig:octagonsurface}
\end{center}
\end{figure}

\begin{example}
Consider the triangular billiard $P$ with boundary given by the triangle with interior angles  $(3\pi/8,\pi/8,\pi/2)$.  We reflect $P$ in the sides emanating from the vertex with angle $\pi/8$. Then, as is expected, we get a surface with 16 copies of the billiard $P$; indeed, $N = \text{lcm}\{8,8,2\}=8$ and so $2N = 16$.  Therefore, the associated generalized polygon is an octagon. Two sides are equivalent if they are translates of each other.  Modding out by this equivalence relation, we have obtained the associated flat surface $\mathcal{S}_P$ as shown in Fig.~\ref{fig:octagonsurface}.  Note, however, that since its genus is greater than $1$, this surface cannot be properly embedded in the plane, but that for convenience, we represent $\mathcal{S}_{P}$ as though it could be done; see Fig.~\ref{fig:octagonsurface}.
\label{exm:triangle}
\end{example}

The surface $\mathcal{S}_P$ is often referred to as a `translation surface' or a `flat surface' in the literature; see, e.g., \cite{Ve3}, \cite{HuSc}, \cite{MaTa}, \cite{Vo} and \cite{Zo}.  Indeed, the geodesic flow on $\mathcal{S}_P$ is nothing but the straight-line flow in the associated generalized polygon, with pairs of opposite sides identified as indicated above.  The key fact concerning the geodesic flow on $\mathcal{S}_P$ is that it is (dynamically) equivalent to the billiard flow on $P$.

\subsection{Singularities of the billiard flow on $P$ and of the geodesic flow on $\mathbf{\mathcal{S}_P}$.} 
\label{subsec:singularities}

As mentioned before, the billiard map is only well defined at certain types of vertices of the rational polygonal billiard. We demonstrate under which conditions the billiard flow can be defined at vertices. If $P$ is a polygon, then the corresponding surface contains copies of the vertices of $P$.   These copies of the vertices are then considered singularities of the flow on the surface and are called \textit{conic singularities}.  

In order to understand what a conic singularity is, we must first introduce the notion of a conic angle.  A \textit{conic angle} is an angle that measures the radians required to form a closed loop about the origin.  In the plane, the conic angle is $2\pi$.  We can form spaces in which the conic angle is not $2\pi$.  Suppose we had a space in which the negative upper half plane was glued to the negative lower half plane and the positive lower  half plane was glued to the positive upper half plane of a second copy of $\mathbb{R}^2$.  Suppose we make similar identifications but instead glue the positive lower half plane to the positive upper half plane of the original copy of $\mathbb{R}^2$.  Then we have formed a space in which the conic angle about the origin is $4\pi$.  Consequently, there are two types of conic singularities: removable and nonremovable.   

\begin{figure}
\begin{center}
\includegraphics{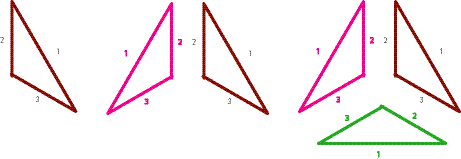}
\caption{The vertex with angle measuring $2\pi/3$ in the triangle above has a conic angle of $4\pi$.  Consequently, the only way to recover a correctly oriented triangle by way of reflecting through sides adjacent to this vertex is to sweep out an angle of $4\pi$.}
\label{fig:weirdgeometry}
\end{center}
\end{figure}

A conic singularity is called \textit{removable} if its  conic angle is $2\pi$.  The conic singularity is called \textit{nonremovable} if its conic angle is $2m\pi$, for some integer $m\geq 2$.   In the context of rational polygonal billiards, we determine the conic singularity of a vertex as follows.  Geometrically, if $v$ is a vertex of a rational billiard corresponding to a nonremovable singularity, when we sweep out an angle of $2\pi$ by continually reflecting the rational billiard in a side, the orientation is not preserved.  In Fig.~\ref{fig:weirdgeometry}, we see that a vertex with an associated conic angle of $4\pi$ must sweep out an angle of $4\pi$ before returning to a copy with the same orientation. We calculate the conic angle of a conic singularity as follows.  If, for $j=1,2,...,r$, $\pi(m_j /n_j)$ is an interior angle of the rational billiard $P$, as in Proposition~\ref{prop:altDefinitionOfRatPolygon}, then the corresponding conic angle is $2m_j\pi$.

Consequently, the geodesic flow on $\mathcal{S}_P$ does not see removable singularities and continues on unimpeded, but is immediately deterred from progressing across the surface when the geodesic intersects a nonremovable singularity.  The geodesic flow on $\mathcal{S}_P$ cannot be determined at such a singularity, implying that the billiard flow on $P$ at the associated vertex cannot be determined either.\footnote{This implication follows from the equivalence of the billiard flow on $P$ and the geodesic flow on $\mathcal{S}_P$.}

\begin{figure}
\begin{center}
\begin{tabular}{ccc}
\includegraphics{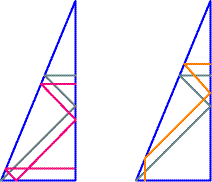}& $\,\,\,\,\,\,$  &\includegraphics{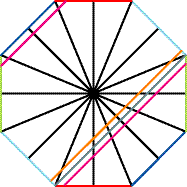}
\end{tabular}
\caption{The surface $\mathcal{S}_P$ corresponding to the triangle $P(3\pi/8,\pi/8,\pi/2)$; the corresponding generalized polygon is an octagon, shown on the right.  Again, we commit the sin of suggesting that the surface can be embedded in the plane, this time to illustrate the ambiguity experienced by the billiard ball at nonremovable singularities.  Once the opposite sides of the octagon have been pairwise identified, as shown in the figure on the right, one obtains the surface $\mathcal{S}_P$.  Note that periodic billiard orbits in the triangle (left) correspond to straight-line paths on the octagon (right), and hence to closed (i.e., periodic) geodesics on $\mathcal{S}_P$.  This figure demonstrates that the two orbits emanating from points near a vertex corresponding to a nonremovable singularity behave in a radically different way; the two orbits are shown in pink and orange, with the grey orbit representing what one would like to believe is the orbit emanating from the nonremovable singularity.}
\label{fig:goodflowonoctagon}
\end{center}
\end{figure}

\begin{example}
Consider again the triangle $P(3\pi/8,\pi/8,\pi/2)$, as was done in Example~\ref{exm:triangle}.  From the above formula, the conic angle of the conic singularity corresponding to the vertex with angle $3\pi/8$ is $6\pi = 2\cdot 3\pi$.  Consequently, this conic singularity is of the nonremovable type.  This means that the flow on the surface $\mathcal{S}_P$ cannot be determined through the corresponding copy of the vertex on $\mathcal{S}_P$.  Moreover, we cannot logically define reflection at this vertex.  The vertex at the center of our surface, however, is of the removable type; the corresponding conic angle is $2\cdot 1 \pi = 2\pi$; see Fig.~\ref{fig:goodflowonoctagon}.
\end{example}

\subsection{$\mathbf{KS_n}$ as a rational billiard.}
\label{subsec:KSNratbilliard}

It was important that we developed an understanding of conic singularities, because, in some sense,  the Koch snowflake is entirely comprised of singularities. Whether or not `corners' of the Koch snowflake\footnote{That is, vertices of the prefractal polygonal approximations $KS_n$, for any $n\geq 0$. It is worth noting that the union of these vertices for all $n\geq 0$ forms a countable dense subset of $\partial(KS)$.} billiard can be shown to correspond to removable or nonremovable singularities of the flow remains to be determined.

\begin{figure}
\begin{center}
\begin{tabular}{cp{20 pt}c}
\includegraphics{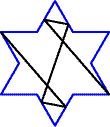} & \vspace{-80 pt} $\Longrightarrow$ & \includegraphics{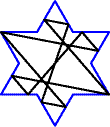}
\end{tabular}
\caption{Acute angles have associated removable conic singularities.  Reflection at acute vertices can then be defined.  We illustrate here a singular orbit of $KS_1$.}
\label{fig:uhoh}
\end{center}
\end{figure}

The prefractal billiard $KS_n$ is a rational billiard.  Indeed, an interior angle of $KS_n$ is either equal to $\pi/3$ (acute) or $4\pi/3$ (obtuse).\footnote{Furthermore, note that except for $n=0$ (when $KS_0 = \Delta$, the equilateral triangle), we always have both removable and nonremovable singularities (and in equal numbers).}  Consequently, $N=\text{lcm}\{3,3\}=3$ and so the corresponding surface contains $2\cdot N = 6$ copies of $KS_n$. Furthermore, removable singularities of the surface correspond to acute angles of $KS_n$ and nonremovable singularities correspond to obtuse angles of $KS_n$.  This means that when $\alpha$ is a singular orbit of $KS_n$, this orbit can be naturally extended to a periodic orbit if the singularities encountered in $KS_n$ were of the removable type; see Fig.~\ref{fig:uhoh}.  Likewise, if $\alpha$ is a singular orbit and the billiard ball encounters an obtuse angle of $KS_n$, the flow on the billiard $P$ cannot be naturally extended.

\begin{figure}
\begin{center}
\includegraphics{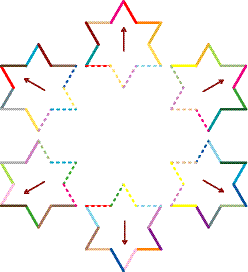}
\end{center}
\caption{The flat surface $\mathcal{S}_{KS_1}$ corresponding to $KS_1$; sides identified. This surface is independent of the flow and is the true surface, with sides properly identified.  We cannot naively reflect in the adjacent sides of a fixed vertex and subsequently identify sides to produce the surface, because $S_{KS_n}$ has genus $g>1$ for all $n>0$.  There are six copies of the table $KS_1$ in the associated flat surface, because $2\cdot \text{lcm}\{3,3\} = 6$.  Moreover, these six copies are generated by letting the dihedral group $D_3$ act on $KS_1$ with an orientation designated by a vector that is not parallel to any of the sides; having such an orientation is a necessary condition for producing the correct number of copies of $KS_1$.  In Fig.~\ref{fig:thecorrectkochsurfacewithflow} below, we illustrate the associated straight-line flow on $S_{KS_1}$ for a particular periodic orbit of $KS_1$ induced by the Fagnano orbit $\fagnano_0$ of $KS_0 = \Delta$.}
\label{fig:thecorrectkochsurface}
\end{figure}

\begin{figure}
\begin{center}
\includegraphics{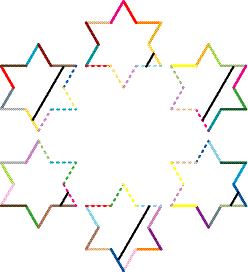}
\end{center}
\caption{The unfolding of the periodic orbit of $KS_1$ induced by the Fagnano orbit $\fagnano_0$ of $KS_0 = \Delta$.}
\label{fig:thecorrectkochsurfacewithflow}
\end{figure}

\begin{figure}
\begin{center}

\includegraphics{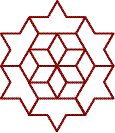} \includegraphics{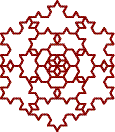} \includegraphics{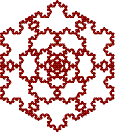}

\end{center}
\caption{The generalized polygons corresponding to $KS_1$, $KS_2$ and $KS_3$, respectively; sides not identified.  These should be viewed as the first, second and third level approximations of the limiting `generalized fractal polygon' presumably associated with the billiard $KS$.  When the proper identifications are made, such surfaces $\mathcal{S}_{KS_1}$, $\mathcal{S}_{KS_2}$ and $\mathcal{S}_{KS_3}$ can be viewed as the first, second and third level approximations of the limiting `fractal flat surface' (as hypothesized in part (i) of Conjecture~\ref{conj:fractalbilliardFractalSurface}).}
\label{fig:threekochsurfaces}
\end{figure}

In Fig.~\ref{fig:thecorrectkochsurface}, we have represented the flat surface $\mathcal{S}_{KS_1}$, and in Fig.~\ref{fig:threekochsurfaces}, we have depicted the generalized polygons corresponding to $\mathcal{S}_{KS_1}$, $\mathcal{S}_{KS_2}$ and $\mathcal{S}_{KS_3}$.  With sides properly identified, the latter generalized polygons become flat surfaces and can be viewed as the first three prefractal approximations of the `fractal flat surface' (of infinite genus) $\mathcal{S}_{KS}$, which we will conjecture (in Section~\ref{sec:conjectures}) to exist as a proper mathematical object; see part (i) of Conjecture~\ref{conj:fractalbilliardFractalSurface}.

\subsection{Equilateral triangle billiards and implications for $\mathbf{KS_n}$.}
\label{subsec:baum}

\begin{figure}
\begin{center}
\includegraphics{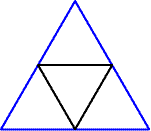}
\caption{The Fagnano orbit $\fagnano_0$ of the equilateral triangle billiard $\Delta = KS_0$.  It is the unique shortest closed billiard trajectory of $\Delta$ and is given by the equilateral triangle inscribed in $\Delta$. The associated initial conditions are $(x_{mid},\pi/3)$, where $x_{mid}$ is the midpoint of one of the sides of $\Delta$ and $\pi/3$ is the corresponding initial angle of the path.}
\label{fig:fagnano0}
\end{center}
\end{figure}

The equilateral triangle billiard $\Delta$ is a very nice rational billiard.  The billiard dynamics on $\Delta$ are well understood.  Two orbits of an equilateral triangle billiard are said to be \textit{equivalent} if each orbit has an unfolding that is a translate of the other and each have the same length.  Denote the equivalence class of an orbit $\alpha$ by $[\alpha]$. In \cite{BaUm}, Baxter and Umble give a classification of the periodic orbits of $\Delta$, up to this notion of equivalence.  Roughly speaking, they partition periodic orbits into two collections, those periodic orbits with odd period and periodic orbits with even period.\footnote{The period of a billiard orbit is the number of times which the pointmass collides with the boundary.  Consequently, an even period orbit is one for which the number of collisions with the boundary is even, and similarly for an odd period orbit.}  Among orbits with even period, they then partition the orbits based on special criteria.  If $\alpha$ is an even periodic orbit, then its equivalence class $[\alpha]$ has the cardinality of the continuum.  The collection of odd period orbits is a collection of odd iterates of the Fagnano orbit $\fagnano_0$, $\{\fagnano_0^{2k+1}\}_{k=0}^\infty$.  Recall that the \textit{Fagnano orbit} $\fagnano_0$ is defined to be the shortest periodic orbit of the equilateral triangle billiard $\Delta$;\footnote{The name ``Fagnano'' for the shortest orbit comes from the name of the Italian mathematician Giovanni Fagnano whose work dates back to 1775.  Indeed, Fagnano researched the existence of the shortest inscribed polygons of a fixed polygon $P$.} it is the unique such closed path and is represented in Fig.~\ref{fig:fagnano0}.    Specifically, for every $k\geq 0$, there is only one element of the equivalence class $[\fagnano_0^{2k+1}]$, that element being $\fagnano_0^{2k+1}$, the ($2k+1$)th repetition of $\fagnano_0$.

\begin{figure}
\begin{center}
\includegraphics{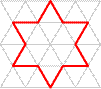}
\end{center}
\caption{Consider a tiling of the plane by the equilateral triangle $\Delta_n$ with side lengths $1/3^n$.  Denote this tiling by $T_{\Delta_n}$. Then we can see that $KS_n\subseteq T_{\Delta_n}$.  This observation is illustrated here for the case when $n=1$.}
\label{fig:singleembedding}
\end{figure}

\begin{figure}
\begin{center}
\includegraphics{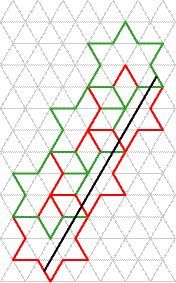}
\end{center}
\caption{As noted in the caption of Fig.~\ref{fig:singleembedding}, $KS_n\subseteq T_{\Delta_n}$. If we consider a periodic orbit $\gamma_{n,i}$ in $\gamma_n$ as shown, then we see that $\gamma_{n,i}$ unfolds to a straight line.  This unfolding is collinear with an unfolding of the Fagnano orbit $\fagnano_0$ of $\Delta_n$.  (See \S\ref{subsec:gammanorbits} for the definition of $\gamma_n$.) As in Fig.~\ref{fig:singleembedding}, the situation is illustrated here for the case when $n=1$.}
\label{fig:singleembeddingwithflow}
\end{figure}

Roughly speaking, the billiard $KS_n$ is comprised of many copies of $\Delta_n$, an equilateral triangle billiard with side lengths measuring $1/3^n$.  More precisely, if $T_{\Delta_n}$ is a tiling of the plane by equilateral triangles with side lengths measuring $1/3^n$, then $KS_n$ can be embedded in $T_{\Delta_n}$; see Fig.~\ref{fig:singleembedding}.  Because of this observation and the fact that all orbits of $\Delta$ can be unfolded in $T_{\Delta_n}$ (see Fig.~\ref{fig:singleembeddingwithflow}), it is reasonable to expect that periodic and singular orbits of $\Delta$ extend naturally to $KS_n$.  In fact, Fig.~\ref{fig:singleembeddingwithflow} illustrates that an unfolded periodic orbit of $KS_1$ is collinear with the unfolding of the Fagnano orbit $\fagnano_0$ of $\Delta=KS_0$.

\section{Experimental Results}
\label{sec:expresults}

We want to describe the `anatomy' of the prefractal billiard $KS_n$.  For fixed $n\geq 0$, $KS_n$ is just a rational polygon.  However, as $n$ changes, the boundary of our prefractal billiard $KS_n$ changes.  We want to be able to talk more concretely about the parts of $KS_n$ that `sprout up' and those that `die off' as $n$ increases.  Parts of $KS_n$ that `sprout up' are what we will call \textit{cells} of $KS_n$.  The collection of segments of $KS_{n-1}$ that are removed to construct the cells of $KS_n$ will then be called the \textit{ghost} of $KS_{n-1}$.  The ghost of $KS_{n-1}$ will not  be present in $KS_n$, but such language will help us discuss the behavior of a certain type of orbits.  

\begin{figure}
\begin{tabular}{p{3.5 cm}  p{3.5 cm}  p{3.5 cm}  p{4 cm}}

\begin{center}\includegraphics{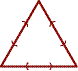}\end{center} & \begin{center}\includegraphics{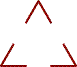}\end{center} & \begin{center}\includegraphics{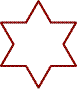}\end{center} & 
\begin{center}\includegraphics{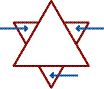}\end{center} \\
(a) The ghost set of $KS_0 = \Delta$, denoted by $G_0$. &  
(b) The elements of the ghost set $G_0$ are removed. &  
(c) Out of every side there `sprouts' two segments, giving rise to $KS_1$.  &  
(d) $G_0\cup KS_1$.  The blue arrows indicate the cells $C_{1,i}$, $1\leq i \leq 3\cdot 4^0=3$, of $KS_1$.\\
\\ 
\end{tabular}
\caption{An illustration of Definitions \ref{def:anatomyOfKSn}(i)--(iii) in terms of $KS_0 = \Delta$ and $KS_1$.  The ghost of the segment $s_{0,i}$, denoted by $g_{0,i}$, is a middle-third segment of $s_{0,i}$ and is removed from $s_{0,i}$ so that we may construct the cell $C_{1,i}$ of $KS_{1}$.  Then $g_{0,i}$ is referred to as the ghost of the cell $C_{1,i}$ of $KS_{1}$.}
\label{fig:illustrateDefis}
\end{figure}

We want to stress that we will not be discussing orbits of $KS_n$.  Rather, we will be discussing what we call \textit{induced orbits} of $KS_n$.  What we mean by `induced' is that if $(x_0,\theta_0)$ is an initial condition of an orbit of $\Delta$, an induced initial condition will be $(x_0^\prime,\theta_0)$, where $x_0$ and $x_0^\prime$ are collinear in the direction of $\theta_0$ and $x_0^\prime$ is a point on $KS_n$.  The new initial condition $(x_0^\prime, \theta_0)$ is what we call an induced orbit of $KS_n$.  The reason for such language is that an orbit of $KS_n$ eventually induces an orbit of $KS_{n+k}$ for some $k\geq 1$.  Such an orbit is induced when the billiard ball passes through a `ghost' of $KS_{n+k-1}$ into a cell of $KS_{n+k}$.  We will discuss this in detail in \S\ref{subsec:gammanorbits}, but first we formally define what we have discussed thus far and give an illustration of these concepts in Fig.~\ref{fig:illustrateDefis}.

\begin{definition}[The anatomy of $KS_n$] 
\label{def:anatomyOfKSn}
{\ }
\renewcommand{\labelenumi}{(\roman{enumi})}
\begin{enumerate}
\item{(Ghosts of $KS_n$). Let $n\geq 0$ and $\{s_{n,i}\}_{i=1}^{3\cdot 4^n}$ be the collection of segments comprising the boundary $\partial(KS_n)$ of the billiard $KS_n$.  Then, for $1\leq i \leq 3\cdot4^n$, the open  middle third of the segment $s_{n,i}$ is denoted by $g_{n,i}$ and is called the \textit{ghost of the segment} $s_{n,i}$.  Moreover, the collection $G_n = \{g_{n,i}\}_{i=1}^{3\cdot 4^n}$ is called the \textit{ghost set of} $KS_n$.  The segments $g_{n,i}$ are removed in order to generate $KS_{n+1}$; see Figs.~\ref{fig:illustrateDefis}(a)--(c).}

\item{(Cell of $KS_n$). Let $n\geq 1$, $1\leq i\leq 3\cdot 4^{n-1}$ and $G_{n-1}$ be the ghost set of $KS_{n-1}$.  Consider the set $G_{n-1}\cup KS_n$ and the bounded region in $G_{n-1}\cup KS_n$ given by an equilateral triangle with side lengths measuring $1/3^n$.  Then this bounded region is what we call a \textit{cell of} $KS_n$.  We denote a cell of $KS_n$ by $C_{n,i}$; see Fig.~\ref{fig:illustrateDefis}(d) (by definition, $KS_0=\Delta$ has no cell).}

\item{(Ghost of a cell). Let $n\geq 1$, $1\leq i\leq 3\cdot 4^{n-1}$ and $C_{n,i}$ be a cell of $KS_n$. Then the ghost $g_{{n-1},i}$ of the segment $s_{n-1,i}$ is called the \textit{ghost of the cell} $C_{n,i}$.  In other words, the ghost of the cell $C_{n,i}$ is the middle-third segment $g_{n-1,i}$ of the segment $s_{n-1,i}$ of $\partial{(KS_{n-1})}$ that is removed from $s_{n-1,i}$ in order to generate that portion of $\partial{(KS_n)}$; see the caption of Fig.~\ref{fig:illustrateDefis}.}
\end{enumerate}
\end{definition}




\begin{example}
If we consider the base of the equilateral triangle with sides having unit length, then the ghost of the base $s_{0,1}$ of $\Delta = KS_0$ is given by $g_{0,1}=(1/3,2/3)$; see Fig.~\ref{fig:illustrateDefis}(a).  
\end{example}

\begin{definition}[Compatible sequence of orbits]
\label{def:compatiblesequence}
For each $n\geq 0$, let $\alpha_n$ be a periodic orbit of $KS_n$.  Then, the sequence $\{\alpha_n\}_{n=0}^\infty$ is said to be \textit{compatible} if each $\alpha_n$ is induced by $\alpha_0$, the initial orbit of $\Delta$.
\end{definition}




We performed computer simulations of induced orbits of the billiard $KS_n$ for $n=1,2,3$.  As a result, we now proceed to pose, discuss and sometimes answer various research questions regarding the induced orbits of $KS_n$, for $n\geq 0$.

\subsection{$\mathbf{\gamma_n}$ Orbits of $\mathbf{KS_n}$.}
\label{subsec:gammanorbits}

\begin{figure}
\begin{center}
\includegraphics{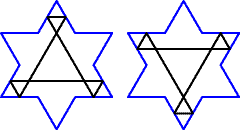}
\end{center}
\caption{There are two orbits in $\gamma_1$.}
\label{fig:KS1Fagnanos}
\end{figure}

We define $\gamma_n$ to be the collection of periodic orbits of $KS_n$ with an initial condition of the form $(x_{mid},\pi/3)$, where $x_{mid}$ is a midpoint of an arbitrary side of $KS_n$.  When $\pi/3$ is not a feasible direction, let $\pi/3$ be an angle relative to the initial side.  Every orbit in $\gamma_n$ can then be demonstrated to be an orbit induced by some periodic orbit of the equilateral triangle $\Delta$; see Fig.~\ref{fig:singleembeddingwithflow} for the case when $n=1$, and for an illustration of the orbits in $\gamma_1$, see Fig.~\ref{fig:KS1Fagnanos}.

We want to know what constitutes a shortest orbit of $KS_n$.  Is the shortest orbit an element of $\gamma_n$?  If not, are all the orbits in $\gamma_n$ of the same length?  Can we justify calling an orbit ``Fagnano'' in some sense that keeps with the `shortest length' meaning of the name?  We now proceed to answer these questions by examining our experimental results.  Out of this will come a clear definition of Fagnano orbit of $KS_n$ and a better understanding of the behavior of the orbits in $\gamma_n$.  

\begin{figure}
\begin{center}
\includegraphics{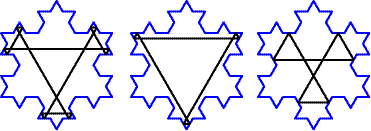}
\end{center}
\caption{We see that as we construct $KS_2$ from $KS_1$, the orbit of $KS_1$ induces an orbit of $KS_2$.  Then, constructing $KS_3$ from $KS_2$, the induced orbit induces an orbit of $KS_3$.  However, this need not always be the case. Certain orbits would remain unchanged as we construct $KS_{n+1}$ from $KS_n$ for some $n$.  Moreover, we can see that the trajectory in the first image passes through the midpoint of the ghost $g_{k,i}$ of a cell $C_{k+1,i}$ of $KS_{k+1}$ for $k \leq 1$.  This fact will be significant in \S\ref{sec:conjectures}.}
\label{fig:KS2Fagnanobreakingup}
\end{figure}

Let $\gamma_{1,i}\in \gamma_1$.  Consider the unfolding of $\gamma_{1,i}$ induced by the classic Fagnano orbit $\fagnano_0$ of $KS_0=\Delta$, and consider the orbit $\alpha$ induced by the initial condition $(x_{mid},\pi/6)$.  Then, the orbit $\alpha$ is clearly shorter than the orbit $\gamma_{1,i}$.  Now, let us also consider the orbit $\gamma_{2,j}\in\gamma_2$ induced by $\fagnano_0$.  Examine the unfolding of $\gamma_{1,i}$ and the unfolding of an equivalent orbit that does not reside on the ghost of the initial side, but on some other one-third-segment of the initial side.  Since equivalent orbits have the same lengths, we can see that this is also a periodic orbit of $\Delta$ that is shorter than $\gamma_{2,j}$ but still longer than $\alpha$.  From this, we deduce that no orbit in $\gamma_2$ will qualify as a shortest orbit of $KS_2$.  Consequently, for every $n>0$, no element of $\gamma_n$ is a shortest orbit of $KS_n$.  However, this does not preclude us from identifying the shortest orbits in $\gamma_n$; see Fig.~\ref{fig:KS2Fagnanobreakingup}.

\begin{figure}
\begin{center}
\includegraphics{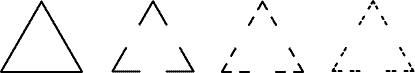}
\end{center}
\caption{From left to right: $\Delta \cap KS_0$, $\Delta \cap KS_1$, $\Delta \cap KS_2$, $\Delta \cap KS_3$. For each approximation $KS_k$ with $k\leq 3$, we are removing a subset of the ghost set $G_k$ of $KS_k$.  The set $\Delta\cap KS$ is comprised of three ternary Cantor sets.} 
\label{fig:intersection}
\end{figure}

Consider the subcollection of $\gamma_n$ comprised of orbits with initial point $x_{mid}\in \Delta\cap KS_n$; see Fig.~\ref{fig:intersection}.  In general, our results indicate that this subcollection is the collection of orbits with shortest length among all orbits in $\gamma_n$. We denote this subcollection by $\fagnano_n$ and call it the collection of Fagnano orbits of $KS_n$.  Recall from our earlier discussion in \S\ref{subsec:baum} that $\fagnano_0$ consists of a single element, namely, the Fagnano orbit (also denoted $\fagnano_0$); see Fig.~\ref{fig:fagnano0}.

We give a special name to the periodic orbit induced by $\fagnano_0$ in the prefractal billiard $KS_n$ (for each fixed $n\geq 0$).  We call this particular orbit the \textit{primary piecewise Fagnano orbit} of $KS_n$ (denoted by $pp\fagnano_n$) because, for every $k\leq n$, the initial segment enters a cell $C_{k,i}$ of $KS_k$ (by passing through the ghost $g_{k-1,i}$ of the cell $C_{k,i}$, as defined in Definition~\ref{def:anatomyOfKSn}(iii)) and, upon reflecting, subsequently forms a true Fagnano orbit in a cell $C_k$.  In general, based on this terminology, we refer to an orbit in $\gamma_n$ as a \textit{piecewise Fagnano orbit}.  In Section 4, we will often denote the collection of piecewise Fagnano orbits of $KS_n$ by  $p\fagnano_n$ (instead of $\gamma_n$).

\subsection{Nonsingular non-$\mathbf{\gamma_n}$ periodic orbits of  $KS_n$.}
\label{subsec:nongammanorbits}

\begin{figure}
\begin{center}
\begin{tabular}{cc}
\includegraphics{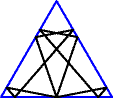}&\includegraphics{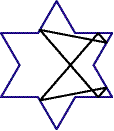}\\
\includegraphics{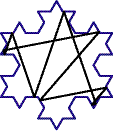}&\includegraphics{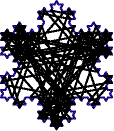}
\end{tabular}
\caption{In this figure, we see an example of a non-$\gamma_n$ periodic orbit of $KS_0$, $KS_1$, $KS_2$ and $KS_3$, where we consider the same inducing condition $(x_{mid},\theta_0)$ in each iteration of the prefractal billiards.  For a fixed inducing condition $(x_0,\theta_0)$, other simulations (not shown) are suggesting that induced orbits of $KS_n$ are (qualitatively) quite convoluted.}
\label{fig:nongammaN}
\end{center}
\end{figure}

Our characterization of periodic orbits of $KS_n$ in terms of what we have referred to as ``$\gamma_n$'' and ``non-$\gamma_n$'' orbits is by no means a rigorous classification of the induced orbits of $KS_n$.  Our results indicate that elements in $\gamma_n$ are not always induced by the Fagnano orbit of $\Delta$ and non-$\gamma_0$ orbits of $\Delta$ do not always induce elements in $\gamma_n$; see Fig.~\ref{fig:nongammaN}.  Since non-$\gamma_0$ orbits of $\Delta$ can induce $\gamma_n$ orbits of $KS_n$, is it possible that non-$\gamma_0$ orbits of $\Delta$ also induce singular orbits of $KS_n$?  Can we demonstrate that there is some nonsingular non-$\gamma_n$ orbit that induces a nonsingular periodic orbit of $KS_{n+k}$, for all $k\geq 1$?

\begin{figure}
\begin{center}
\begin{tabular}{cc}
\includegraphics{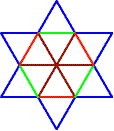}\\
\includegraphics{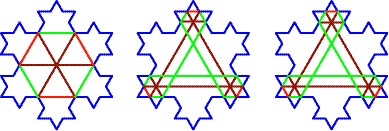}
\end{tabular}
\caption{The billiard flow cannot be determined at obtuse angles.  However, examination of the corresponding geodesic flow on $\mathcal{S}_{KS_n}$ for $n=1,2$ allows us to draw  `possible' paths; note the ambiguity at obtuse angles.  It is intriguing that the `possible' choices that the billiard ball `makes' at the triadic points $\{1/3,2/3\}$ and $\{1/9,2/9\}$ are exactly the same.  This observation may allow us to force a logically consistent definition of reflection at nonremovable singularities of the billiard map on $KS_n$.  In other words, this figure suggests that there may be a rigorous way of reducing the ambiguity the billiard ball experiences at obtuse angles of $KS_n$.  Because of the unique `infinite' symmetry enjoyed by the snowflake boundary $\partial (KS)$, we may be able to express the billiard flow on $KS$ at nonremovable singularities in terms of a (countably) infinite number of `possible' paths the point mass may take.   Again, we may be aided in this investigation by considering the (yet to be determined) geodesic flow on the associated `fractal flat surface' $\mathcal{S}_{KS}$; see Conjecture~\ref{conj:fractalbilliardFractalSurface}.}
\label{fig:gammanAtObtuseAngles}
\end{center}
\end{figure}

\begin{figure}
\begin{center}
\begin{tabular}{cc}
\includegraphics{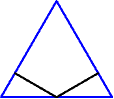}&\includegraphics{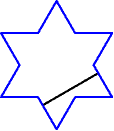}\\
\includegraphics{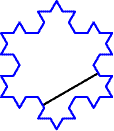}&\includegraphics{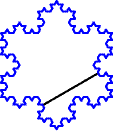}
\end{tabular}
\caption{In this figure, we see an example of a non-$\gamma_n$ periodic orbit of $KS_0$, $KS_1$, $KS_2$ and $KS_3$.  We consider the same inducing condition $(x_{mid},\pi/6)$ in each iteration of the prefractal billiards.}
\label{fig:KSpi6}
\end{center}
\end{figure}

If $n=0$ and $x_0\neq x_{mid}$, then $\alpha_{(x_0,\pi/3)}$, the periodic orbit of $KS_n$ with initial (or `inducing') condition $(x_0,\pi/3)$, is certainly a permissible element of the collection of non-$\gamma_n$ orbits of $KS_n$.  If, for some integers $j,\nu\geq 1$ such that $0< \nu < 3^j$,  we let $x_0=(\nu/3^j,0)$ on the base of $\Delta$, then $\alpha_{(x_0,\pi/3)}$ induces a singular orbit of $KS_n$ for all $n\geq j$; see Fig.~\ref{fig:gammanAtObtuseAngles}.  Moreover, this singular orbit collides with an obtuse angle of $KS_n$, meaning that we cannot make sense of the billiard flow with this particular initial condition.  However, experimental results are indicating that for an initial condition $(x_0,\theta_0) = (x_{mid},\pi/6)$, the resulting orbit (i) will be periodic and (ii) may be nonsingular, for all $n\geq 0$; see Fig.~\ref{fig:KSpi6}.

\begin{figure}
\begin{center}
\begin{tabular}{cc}
\includegraphics{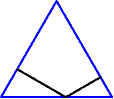}&\includegraphics{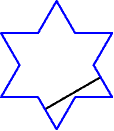}\\
\includegraphics{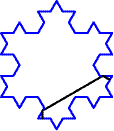}&\includegraphics{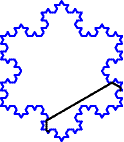}\\
\includegraphics{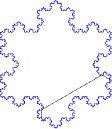}&\includegraphics{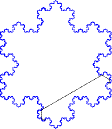}\\
\end{tabular}
\caption{In this figure, we see an example of what happens when we shift the initial point from $x_{mid}$ to $x_{mid}+\delta x$, for some small perturbation $\delta x$.   We consider the same inducing condition $(x_{mid}+\delta{x},\pi/6)$ in each iteration of the prefractal billiards.  What is interesting in this particular example is that the billiard trajectory seems to want to enter into a new cell with each iteration of the prefractal billiard.  This indicates that there may be a limiting object to consider, presumably, a billiard trajectory of $KS$.}
\label{fig:KSpi6shifted}
\end{center}
\end{figure}

With regards to the billiard $\Delta$, changing the basepoint of the initial condition $(x_0,\theta_0)$ will not affect the periodicity of the orbit, but may affect the length of the orbit.  In our experiments regarding orbits in $\gamma_n$, we observed that changing the basepoint did not affect the periodicity of the orbit or the length of the orbit.  This can be seen by unfolding the orbits.  However, our experiments on non-$\gamma_n$ orbits indicated otherwise. One would like to see that orbits emanating from a side $s_{n,i}$ with the same direction are equivalent, but this did not turn out to be the case.  In particular, an example of a nonsingular non-$\gamma_n$ periodic orbit that increased in length after a translation of the basepoint is one given by the inducing condition $(x_{mid},\pi/6)$; see Fig.~\ref{fig:KSpi6}.  Unlike the case of $\gamma_n$-type orbits, if we shift the basepoint of the inducing condition by some small perturbation $\delta x$, then the resulting orbit is longer in length; see Fig.~\ref{fig:KSpi6shifted}.

\subsection{Singular orbits of $KS_n$.}
\label{subsec:singularorbits}
Recall that a singular orbit of $\Delta$ is an orbit which collides with a vertex of the boundary $\partial \Delta$.  As we discussed in \S\ref{subsec:singularities}, vertices of a billiard with associated conic angles measuring $2\pi$ radians are considered removable singularities of the geodesic flow on $\mathcal{S}_P$.  Since every vertex in $\Delta$ has an associated removable conic singularity, the billiard flow may be naturally defined at corners of $\Delta$.  

Let $\alpha$ be an induced singular orbit of $KS_n$.  Since $\alpha$ is a singular orbit of $KS_n$, $\alpha$  collides with either an obtuse angle or an acute angle of $KS_n$.  Recall that acute angles of $KS_n$ constitute removable singularities of the billiard flow whereas obtuse angles of $KS_n$ constitute nonremovable singularities of that flow. Consequently, we may not be able to always determine how $\alpha$ behaves at vertices of $KS_n$.  As $n$ increases, $KS_n$ grows. Therefore, we may ask whether or not an induced singular orbit $\alpha$ remains singular as $n$ increases.  Because $\alpha$ can only be naturally extended at an acute angle vertex of $KS_n$, is it possible to prove that induced singular orbits always collide with acute angles? Moreover, do singular orbits which only collide with acute angles at some stage $n$ do so at every later stage $n+k$, $k\geq 1$?  Even though our attention is focused on induced orbits, is it possible that there exist singular orbits of $KS_n$ that are not induced by any orbit of $\Delta$?	

\begin{figure}
\begin{center}
\begin{tabular}{cc}
\includegraphics{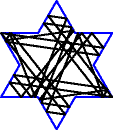}&\includegraphics{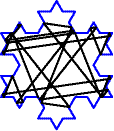}
\end{tabular}
\end{center}
\caption{The induced singular orbit only collides with acute angles of $KS_1$.  The same initial condition induces a singular orbit of $KS_2$, whereas this time the orbit collides with obtuse angles of $KS_2$. Insufficient numerical accuracy prevents us from investigating (via computer simulations) the induced singular orbit of $KS_3$.}
\label{fig:mixedvertices}
\end{figure}

If $\alpha$ is an induced singular orbit of $KS_n$, then experimental results indicate that $\alpha$ will always have to collide with a vertex of $KS_{n+k}$, for all $k\geq 0$.  An explanation for why this may be so is that $KS_{n+k}$ can be embedded in $T_{\Delta_{n+k}}$ and the inducing initial condition gave rise to a singular orbit of $\Delta$.  However, our experimental results show that there is no consistency in which types of vertices are encountered by the singular orbit.  In particular, it can happen that a singular orbit of $KS_n$ may only collide with acute angles, but the induced orbit of $KS_{n+k}$ will collide with an obtuse angle for some $k\geq 1$ and thus end there; see Fig.~\ref{fig:mixedvertices} for a concrete example.

We note that because $KS_n$ can be embedded in $T_{\Delta_n}$ (see Fig.~\ref{fig:singleembedding}), it is reasonable to expect that every singular orbit of $KS_n$ is induced.

\subsection{Quasiperiodic orbits of $\mathbf{KS_n}$.}
\label{subsec:quasiperiodic}

When one considers a rational billiard $P$, the natural question to ask is ``\textit{What are the periodic orbits of} $P$?'' One usually does not consider initial conditions that do not give rise to periodic orbits.\footnote{See, however, the well-known Veech dichotomy [\textbf{Ve1--3}] (as described, e.g., in [\textbf{HuSc},\textbf{MaTa},\textbf{Vo}]), according to which a direction is either periodic or  uniquely ergodic in the rational billiard table.} Moreover, if one has identified periodic orbits of $P$, aperiodic ``approximations'' to these orbits may not be worth analyzing.  Nevertheless, when the boundary is changing---as is the case in the recursive construction of the Koch snowflake billiard table $KS$ via the rational polygonal billiard tables $KS_n$---such  approximations, however such a notion is defined, may be interesting. We call these approximate orbits \textit{quasiperiodic orbits}.  Recall that every orbit of a billiard is assumed to have a unit-speed parameterization.  Formally, we define what we mean by a quasiperiodic orbit of $KS_n$ as follows.

\begin{definition}[Quasiperiodic orbits of $KS_n$]
Fix $T>0$, $n\in \mathbb{N}$, and let $\mathscr{P}_n$ denote the collection of all periodic orbits of $KS_n$.  An orbit $q$ of the billiard $KS_n$ is called \textit{quasiperiodic} if given any $\epsilon >0$, there exists a finite partition of $[0,T]$, $t_0=0 < t_1 < t_2 < ...< t_l=T$, and periodic orbits $p_j\in \mathscr{P}_n$ such that
\begin{eqnarray*}
	|q(t)-p_j(t)|<\epsilon, \ \ \text{for all} \ \ t\in [t_{j-1},t_{j}] \ \ \text{and} \ \  j\in \{1,2,...,l\}.
\end{eqnarray*}
\label{def:quasiperiodicorbit}
\end{definition}

\begin{figure}
\begin{center}
\begin{tabular}{ccc}
\includegraphics{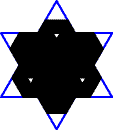} & \includegraphics{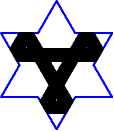} & \includegraphics{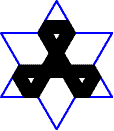}
\end{tabular}
\end{center}
\caption{The left image depicts the quasiperiodic orbit $q$.  The next two images depict the behavior of $q$ over two subintervals of time $(t_0,t_1)$ and $(t_1,t_2)$.  In the middle image, the corresponding periodic orbit $p_1$ would be the one induced by the initial condition $(x_{mid},\pi/3)$ and the corresponding periodic orbit $p_2$ would be induced by $((1/4,0),\pi/3)$, with the base of the generating equilateral triangle $\Delta$ lying on the $x$-axis and the lower left vertex being located at the origin.}
\label{fig:quasiperiodic}
\end{figure}

As the definition would indicate, quasiperiodic orbits exist in any rational polygonal billiard table, but may have special significance in the Koch snowflake billiard.  In Fig.~\ref{fig:quasiperiodic}, we see an example of a quasiperiodic orbit of $KS_1$ that is `shadowing' the $\gamma_1$ orbits of $KS_1$.  In general, if we consider a rational approximation to the initial angle $\theta_0$, we obtain a  quasiperiodic orbit of $KS_n$. If $n$ is changing and the billiard ball is passing through the ghost of a cell $C_{n+k,i}$, then we want to have a way of discussing the `stability' of the resulting induced quasiperiodic orbits of $KS_{n+k}$, for all $k\geq 1$.

\begin{definition}[Stability of a sequence of quasiperiodic orbits]
\label{def:stablequasiperiodicorbit}  
Let $q_0$ be a quasiperiodic orbit of $\Delta$ and $p_0$ a periodic orbit shadowed by $q_0$ over some subinterval $[t_{j-1},t_{j}]\subseteq [0,T]$ (as in Definition~\ref{def:quasiperiodicorbit}).  Then we say that the sequence $\{q_n\}_{n=0}^\infty$ of quasiperiodic orbits of $KS_n$ induced by $q_0$ is \textit{stable} if the sequence $\{p_n\}_{n=0}^\infty$ of periodic orbits induced by $p_0$ has the property that $q_n$ shadows $p_n$ for every $n\geq 0$.\footnote{Hence, in the sense of Definition~\ref{def:compatiblesequence} above, $\{q_n\}_{n=0}^\infty$ is assumed to be a `compatible' sequence of quasiperiodic orbits.}
\end{definition}

With regards to quasiperiodic orbits of $KS_n$, there are two types of potential behaviors that we would like to investigate.  Specifically, we are naturally led to ask the following two questions.  Firstly, if we fix $n\geq 0$ as well as a given periodic orbit $p$ of $KS_n$, is there a sequence of quasiperiodic orbits $\{q_i\}_{i=1}^\infty$  of $KS_n$ that converges to $p$, and if so, in what sense?  Secondly, for a fixed quasiperiodic orbit $q$ of $KS_n$, in what ways is the induced quasiperiodic orbit of $KS_{n+1}$ qualitatively and quantitatively different from $q$?  More precisely, if $q_n$ is a quasiperiodic orbit of $KS_n$, $n\geq 0$, with a fixed rational angle $a/b$ independent of $n$, and $\{q_{n+k}\}_{k=0}^\infty$ is a sequence of induced quasiperiodic orbits, is the latter sequence stable (in the sense of Definition~\ref{def:stablequasiperiodicorbit})?

\begin{figure}
\begin{center}
\includegraphics{quasiperiodic1and2.png} \includegraphics{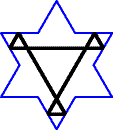} \includegraphics{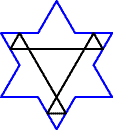}
\caption{This figure depicts what happens as the approximate angle $\theta_0$ becomes more and more accurate.  The quasiperiodic orbit seems to dry up as the initial angle $\theta_0\to\pi/3$.  Here, like elsewhere in the experimental results reported in \S\ref{subsec:quasiperiodic}, the initial angle $\theta_0$ is a rational approximation of the initial angle $\theta$ of a given periodic orbit $p$, obtained via the continued fraction expansion of $\theta$.}
\label{fig:quasiDryingUp}
\end{center}
\end{figure}

\begin{center}
\begin{figure}
\begin{tabular}{cc}
\includegraphics{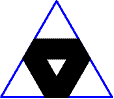} & \includegraphics{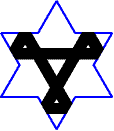} \\
 \includegraphics{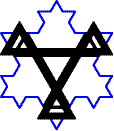} & \includegraphics{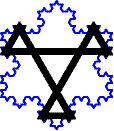}
\end{tabular}
\caption{As the snowflake ``grows'', we see that the quasiperiodic orbit breaks up over the boundary, but does not devolve into an incoherent collection of segments. Rather, there is some sense of stability.  That qualitative sense of stability is captured quantitatively in Definition~\ref{def:stablequasiperiodicorbit}.}
\label{fig:quasiperiodicBreakingUp}
\end{figure}
\end{center}

\begin{center}
\begin{figure}
\begin{tabular}{cc}
\includegraphics{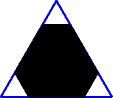} & \includegraphics{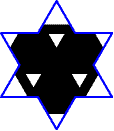} \\ \includegraphics{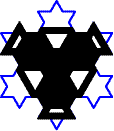} & \includegraphics{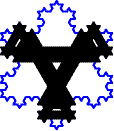}
\end{tabular}
\caption{Comparing this to Fig.~\ref{fig:quasiperiodicBreakingUp}, we see that for fixed $n\geq 0$, a quasiperiodic orbit of $KS_n$ fills more of $KS_n$ as $T$ increases, where $T$ is as stated in Definition~\ref{def:quasiperiodicorbit}.  If we consider a sequence of induced quasiperiodic orbits $\{q_n\}_{n=0}^\infty$, with $q_n$ a quasiperiodic of $KS_n$, then this figure, along with Fig.~\ref{fig:quasiperiodicBreakingUp}, illustrates that the stability of $\{q_n\}_{n=0}^\infty$ (in the sense of Definition~\ref{def:stablequasiperiodicorbit}) is (qualitatively) independent of the time $T$.  Note, however, that both the partitioning of $[0,T]$ into subintervals, and the time intervals over which the shadowing of $p_n$ by $q_n$ occurs, are in general dependent on $n$.}
\label{fig:quasiBreakingUpLong}
\end{figure}
\end{center}

Our experimental results thus far are indicating that we can answer all of these questions in the affirmative.  If $\{a_i/b_i\}_{i=1}^\infty$ is a sequence of rational approximations of $\pi/3$, then Fig.~\ref{fig:quasiDryingUp} describes what happens as $a_i/b_i \to \pi/3$.  This suggests that for a fixed $n$, there is some notion of convergence of quasiperiodic orbits to a corresponding induced periodic orbit $\alpha$ of $KS_n$.  For example, in the present case of Fig.~\ref{fig:quasiDryingUp} (where $n=1$), the corresponding sequence $\{q_i\}_{i=1}^\infty$ of quasiperiodic orbits $KS_1$ seems to be converging (in some suitable sense) to the primary piecewise Fagnano orbit $pp\fagnano_1$ of $KS_1$.  As $n$ increases, Fig.~\ref{fig:quasiperiodicBreakingUp} illustrates that the quasiperiodic orbits break up over the boundary, but that the induced sequence of quasiperiodic orbits appears to remain stable, in the sense of Definition~\ref{def:stablequasiperiodicorbit}.  Furthermore, if we increase the time $T$, as given in Definition~\ref{def:quasiperiodicorbit}, then Fig.~\ref{fig:quasiBreakingUpLong} shows that the same behavior occurs.

\section{Conjectures and Open Problems}
\label{sec:conjectures}

We propose here several conjectures and open problems regarding the billiard $KS$. As we have done in the previous sections, we begin by first defining a few necessary terms.  Our ability to investigate the proposed billiard $KS$ has been limited by the finiteness of the computer.  Hence, our conjectures deal mostly with ``\textit{what happens in the limit}''.  However, it is not exactly clear what we mean by \textit{limit}.  Hence, when we say `limit', we are assuming a suitable notion of limit, most likely (depending on the context) an inverse limit of some inverse system or the Gromov--Hausdorff limit of a sequence of complete metric spaces; see, respectively, \cite{HoYo} or \cite{Gr} for an introduction to these notions. 

Once a suitable notion of limit has been defined, we want to consider the (admittedly, possibly empty) collection of periodic orbits of the billiard $KS$, which we denote by $\mathscr{P}$.  However, our focus thus far  has been on induced orbits of $KS_n$.  We define $\mathscr{I}$ to be the collection of induced periodic orbits of $KS$, and phrase our conjectures in terms of elements in $\mathscr{I}$.  As in Section~\ref{sec:expresults}, where we discussed the `anatomy' of $KS_n$, we now discuss the anatomy of the proposed billiard $KS$ and the collection of induced orbits $\mathscr{I}$.  We define the \textit{footprint}\footnote{This term is sometimes used in the literature on billiards; see, e.g., \cite{BaUm}.} of an orbit to be the points of the orbit that belong to the boundary of the billiard $KS$.  We call the \textit{ghost set} of $KS$ the collection $G=\bigcup_{n=1}^\infty G_n$, the union of all ghost sets $G_n$ associated with $KS_n$ (see Definition~\ref{def:anatomyOfKSn}(i)).

\begin{definition}[Footprint of an orbit]
Let $\alpha$ be an orbit of a billiard $P$.  Then, the \textit{footprint} of $\alpha$ is the set of points of the orbit that belong to the boundary $\partial P$.  
\end{definition}

\begin{definition}[Ghost set of $KS$]
If $G_n$ is the ghost set of $KS_n$, then the \textit{ghost set} of $KS$ is defined by $G=\bigcup_{n=0}^\infty G_n$.
\end{definition}

\begin{definition}[Self-similar orbit]
\label{def:selfsimorbit}
Let $\alpha$ be a periodic orbit of the proposed billiard $KS$.  Then, $\alpha$ is said to be a \textit{self-similar orbit} if its footprint is a self-similar subset of $\partial(KS)\subseteq \mathbb{R}^2$.\footnote{\label{foot:abuseoflanguage}Here and thereafter, when talking about an IFS or a self-similar set, we are making the same slight abuse of language as in \S\ref{sec:introduction} when referring to the Koch snowflake curve $\partial(KS)$ as being `self-similar', whereas strictly speaking, $\partial(KS)$ is obtained by fitting together three copies of a self-similar set, namely, the von Koch curve; see Figs.~\ref{fig:kochcurveconstruction}~and~\ref{fig:3kochcurves}.} 
\end{definition}

Note, in particular, that by construction, the sequence of primary piecewise Fagnano orbits $\{pp\fagnano_n\}_{n=0}^\infty$ is compatible, in the above sense (see Definition \ref{def:compatiblesequence}).

\begin{figure}
\begin{center}
\includegraphics{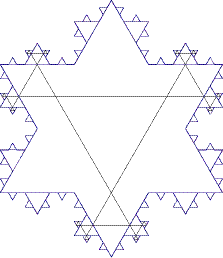}
\end{center}
\caption{Motivation for building the primary Fagnano orbit of $KS_n$ by means of an iterated function system (IFS):  the IFS produces the orbit $pp\fagnano_n$ of $KS_n$ by first contracting, rotating and translating the Fagnano orbit of $\Delta$, so that what results are three copies of $\fagnano_0$ appropriately appended to $\fagnano_0$.  In actuality, there are three IFS's ($F_1,F_2,F_3$) acting in unison to produce this picture. Specifically, $F_1$ acts on the upper left scaled copy of $\fagnano_0$ to produce two scaled, rotated and translated copies of that piece of $pp\fagnano_1$.  Likewise, $F_2$ acts on the upper right copy of $\fagnano_0$ to produce two scaled, rotated and translated copies of that piece of the orbit $pp\fagnano_1$. Finally, $F_3$ behaves similarly on the bottom scaled copy of $\fagnano_0$.  The union of all three images unioned with $pp\fagnano_1$ is then the orbit $pp\fagnano_2$.  Continuing in this manner, we can produce $pp\fagnano_n$, for every $n\geq 0$.}
\label{fig:IFSconstruction}
\end{figure}

Our first conjecture asserts the existence of a very special periodic orbit of the proposed billiard $KS$.  We denote the (presumed) `limit' of the orbits $pp\fagnano_n$ by $pp\fagnano$, and call this the \textit{primary piecewise Fagnano orbit} of $KS$. We can see from Fig.~\ref{fig:IFSconstruction} that there is an iterated function system (IFS),\footnote{More specifically, in the spirit of footnote~\ref{foot:abuseoflanguage}, this `IFS' is given by three different IFS's, denoted by $F_1$, $F_2$, $F_3$ and each comprised of two contractive similarity transformations of $\mathbb{R}^2$; see the caption of Fig.~\ref{fig:IFSconstruction}. The associated attractor is therefore a `self-similar set' of $\mathbb{R}^2$ (also in the spirit of footnote~\ref{foot:abuseoflanguage}).} denoted by $F$, which is producing scaled copies of pieces of the orbit, and such that for any $n\geq 0$, $pp\fagnano_n = \bigcup_{i=0}^n F^i(\fagnano_0)$.  Observe that for each finite $n$, the chronology\footnote{The chronology of an orbit is the order in which the billiard ball visits points of the boundary of the billiard table.} of the orbit $pp\fagnano_n$ can be easily recovered. Furthermore, if we consider the fixed point attractor of the IFS,\footnote{See, e.g., \cite{Fa} for a detailed discussion on iterated function systems.}

\begin{eqnarray*}
\fagnano &=& \lim_{n\to\infty} F^n(\fagnano_0), 
\end{eqnarray*}

\noindent then this set $\fagnano$ should correspond to the footprint of the proposed periodic orbit $pp\fagnano$ of $KS$.\footnote{Note that by construction, $\fagnano$ is a self-similar subset of $\partial(KS) \subseteq \mathbb{R}^2$.}  Defining the orbit $pp\fagnano$ to be $\fagnano$ is not very satisfactory, however, because such a definition does not provide useful information on how to recapture any sense of chronological order.  Accordingly, we make the following formal conjecture.

\begin{conjecture}
\label{conj:weakppf}
Under a suitable notion of limit, we conjecture that the chronology of $pp\fagnano$ can be naturally realized and hence, that $pp\fagnano$ can  be considered as a true periodic orbit of the Koch snowflake billiard $KS$.
\end{conjecture}

Significant analytical evidence in support of $pp\fagnano$ being a well-defined orbit of the billiard $KS$ is the fact that (i) $pp\fagnano$ is a self-similar orbit (specifically, its footprint appears to be a self-similar Cantor set) and (ii) has finite length.  

Once we have established the existence of the periodic orbit $pp\fagnano$, we expect that demonstrating the existence of other periodic orbits given as `limits' of compatible sequences of induced periodic orbits $\{\alpha_n\}_{n=0}^\infty $ may not be an unreasonable endeavor.  Specifically, we conjecture that orbits emanating from the boundary of the proposed billiard $KS$ at an angle of $\pi/3$ can be realized by a suitable generalization of the procedure (yet to be precisely determined) that would substantiate our claim that $pp\fagnano$ forms a periodic orbit of the Koch snowflake billiard.

Our next conjecture makes specific Conjecture~\ref{conj:weakppf} and the comments surrounding it.  It also discusses some of the geometric properties of the presumed `piecewise Fagnano orbits' of $KS$, including the `primary piecewise Fagnano orbit' $pp\fagnano$, the existence of which was asserted in Conjecture~\ref{conj:weakppf}.

\begin{conjecture}[Existence of self-similar periodic orbits of $KS$] 
\label{conj:existenceOfSelfSimilarPeriodicOrbits}
{\ }
\renewcommand{\labelenumi}{(\roman{enumi})}
\begin{enumerate}
\item{\emph{(The primary piecewise Fagnano orbit of $KS$)}. The periodic orbit $pp\fagnano$ of $KS$ can be realized \emph{(}and its chronology restored\emph{)} as a suitable limit \emph{(}possibly, an inverse limit\emph{)} of the compatible sequence of primary piecewise Fagnano orbits $pp\fagnano_n$ of $KS_n$.  Furthermore, the periodic orbit $pp\fagnano$ is a self-similar orbit \emph{(}in the sense of Definition~\ref{def:selfsimorbit}\emph{)}; more specifically, its footprint is the natural self-similar middle-third Cantor set, strung around the boundary $\partial(KS)$.}
\item{\emph{(Piecewise Fagnano orbits of $KS$)}.  More generally, an arbitrary `piecewise Fagnano periodic orbit' $\alpha$ of $KS$ can be defined and realized as follows.  There exists a compatible sequence $\{\alpha_n\}_{n=0}^\infty$ of periodic orbits such that for each $n\geq 0$, $\alpha_n$ belongs to $p\fagnano_n$, and $\{\alpha_n\}_{n=0}^\infty$ converges \emph{(}in a suitable sense\emph{)} to $\alpha$.\footnote{Recall from the end of \S\ref{subsec:gammanorbits} that $p\fagnano_n$ stands for the collection of piecewise Fagnano periodic orbits of $KS_n$.} Furthermore, $\alpha$ is a self-similar orbit of $KS$.}
\end{enumerate}
\end{conjecture}

The following open problem is motivated by our earlier discussion in \S\ref{subsec:nongammanorbits} and \S\ref{subsec:singularorbits}.\footnote{Furthermore, it has an obvious counterpart in the case where $x_0 \neq x_{mid}$ is allowed to be on any of the sides of $\Delta$.}

\begin{openproblem}
On the base of $\Delta = KS_0$, does there exist $x_0$ \emph{(}other than the midpoint $x_{mid}$\emph{)} such that the sequence $\{\alpha_n\}_{n=0}^\infty$ of compatible non-$\gamma_n$ periodic orbits induced by the initial condition $(x_0,\pi/6)$ will \emph{(}i\emph{)} avoid all obtuse angle vertices of $KS_n$, for every $n\geq 0$ and \emph{(}ii\emph{)} converge to a periodic orbit $\alpha$ of $KS$?
\end{openproblem}

Moreover, we conjecture that if  the initial angle $\theta_0$ of the trajectory is not $\pi/3$, but such that $(x_0,\theta_0)$ would have induced a nonsingular non-$\gamma_n$ periodic orbit $\alpha_n$ of $KS_n$ for all $n\geq 0$, then the associated compatible sequence $\{\alpha_n\}_{n=0}^\infty$ `converges' to a periodic orbit $\alpha$  of $KS$; and correspondingly, the associated sequences of footprints of $\{\alpha_n\}_{n=0}^\infty$ converges to the footprint of $\alpha$. In addition, we conjecture that the footprint of $\alpha$ is a topological Cantor set.\footnote{I.e., it is a perfect and totally disconnected subset of the boundary $\partial(KS)$.} A primary candidate for illustrating this conjecture is the periodic orbit induced by the initial condition $(x_{mid} + \delta x, \pi/6)$, for some small perturbation $\delta x$; see Fig.~\ref{fig:KSpi6shifted}. 

A key step towards placing the above conjectures and open problems in a broader context and on firmer mathematical grounds must involve the theory of flat surfaces associated with rational billiards, along with its eventual extension to surfaces of infinite genus (called here `fractal flat surfaces') presumably associated with fractal billiards.\footnote{The flat surface $\mathcal{S}_P$ of a rational polygonal billiard $P$ always has a finite genus. Moreover, a simple calculation based on the known formula for $g(\mathcal{S}_P)$ (see, e.g., Lemma 1.2, p. 1022 of \cite{MaTa}) shows that $g_n = g(\mathcal{S}_{KS_n}) \to\infty$ as $n\to\infty$.  Hence, it is natural to expect that the presumed limiting surface $\mathcal{S}_{KS}$ has infinite genus; see Conjecture~\ref{conj:fractalbilliardFractalSurface}.}    Ultimately, we are conjecturing the existence of a suitably defined billiard table $KS$.  This amounts to establishing that (i) there is a corresponding flat surface $\mathcal{S}_{KS}$ and (ii) the billiard flow on $KS$ is equivalent to the billiard flow on $\mathcal{S}_{KS}$.  In particular, we conjecture that the `limit' of the sequence of prefractal billiard tables $\{KS_n\}_{n=0}^\infty$ would correspond to this billiard $KS$, and that similarly, the associated surface $\mathcal{S}_{KS}$ would be the `limit' of the sequence of  prefractal flat surfaces $\{\mathcal{S}_{KS_n}\}_{n=0}^\infty$.  If we can prove that there are removable singularities in the proposed surface $\mathcal{S}_{KS}$, then it is not unreasonable to expect the conjecture to be true.  Moreover, if the wonderful symmetry of the Koch snowflake curve $\partial(KS)$ can be exploited so as to provide us with a way of dealing with the ambiguity the billiard ball experiences at non-removable singularities (see, e.g., the caption of Fig.~\ref{fig:gammanAtObtuseAngles}), then we would be further justified in believing the conjecture to be a plausible statement.

We summarize the main aspects of this central conjecture---and clearly, long-term open problem---in the following more specific form. Naturally, part of the difficulty in dealing with this problem will involve finding the appropriate notions of limits involved in the formulation of the conjecture.\footnote{For example, even if the notion of `inverse limit' is suitable for formulating part (ii) or the end of part (i) of the conjecture, one would still need to specify the maps involved in the definition of the underlying inverse systems; see, e.g., \cite{HoYo}.}

\begin{conjecture}[The fractal billiard $KS$ and fractal surface $\mathcal{S}_{KS}$, along with the associated flows]
\label{conj:fractalbilliardFractalSurface}
{\ }
\renewcommand{\labelenumi}{(\roman{enumi})}
\begin{enumerate}
\item \emph{(The fractal flat surface $\mathcal{S}_{KS}$, along with the geodesic flow)}. The sequence $\{\mathcal{S}_{KS_n}\}_{n=0}^\infty$ of prefractal flat surfaces associated with the rational billiard $KS_n$ converges \emph{(}in the Gromov--Hausdorff sense, see \cite{Gr}\emph{)} to a surface of infinite genus $\mathcal{S}_{KS}$, called the `fractal flat surface' associated with $KS$.  Correspondingly, the `geodesic flow' on $\mathcal{S}_{KS}$ can be realized as a suitable limit \emph{(}possibly, an inverse limit\emph{)} of the geodesic flows on the surfaces $\mathcal{S}_{KS_n}$. 

\item \emph{(The Koch snowflake fractal billiard $KS$, along with its billiard flow)}. The fractal billiard $KS$ can be defined as a suitable limit \emph{(}possibly, an inverse limit\emph{)} of the prefractal polygonal rational billiards $KS_n$; essentially,\footnote{It is clearly true (and well known) that viewed as a sequence of compact subsets of $\mathbb{R}^2$, the billiard tables $KS_n$ converge in the sense of the Hausdorff metric (or, equivalently, in the Gromov--Hausdorff sense) to the Koch snowflake billiard table $KS$.  Hence, the real issue concerns here the associated billiard flows.} this means that the billiard flow on $KS$ can be realized as a suitable limit \emph{(}also possibly, an inverse limit\emph{)} of the corresponding billiard flows on the billiard tables $KS_n$.
\item \emph{(Geodesic vs. billiard flow)}. Finally, the geodesic flow on the fractal flat surface $\mathcal{S}_{KS}$ is equivalent to \emph{(}and, at first, may be used as a suitable substitute for\emph{)} the billiard flow on the Koch snowflake billiard $KS$.
\end{enumerate}
\end{conjecture}

We close this paper by stating the following very long-term problem, which is directly motivated by the questions raised in \cite{La2} concerning the relationship between `fractal drums' and `fractal billiards', including the Koch snowflake drum and billiard.  See, in particular, \cite{La2}, Conjecture 6, p. 198, itself motivated by Conjectures 2 and 3, pp. 159 and 163--164, respectively; see also \S{12.5.3} of  \cite{La-vF}. For information regarding the Koch snowflake drum and other `fractal drums' (viewed as `drums with fractal boundary'), see, e.g., [\textbf{La1--2},\textbf{LaNRG},\textbf{LaPa}], [\textbf{La-vF},\S{12.3} \& \S{12.5}], and the relevant references therein.  For information regarding trace formulas (including the Gutzwiller and Chazarain trace formulas) in various contexts, we refer, e.g., to \textbf{[Gz1,2]}, \cite{Ch}, \cite{Co} and \cite{DuGn}.

\begin{openproblem}[Fractal Billiard vs. Fractal Drum.]  
\label{openproblem:fractalbilliardvsfractaldrum} Once the existence of the Koch snowflake billiard $KS$ has been firmly established \emph{(}as hypothesized in Conjecture~\ref{conj:fractalbilliardFractalSurface}\emph{)}, can one formulate, and eventually establish, a suitable `fractal trace formula' in this context?  Presumably, the latter would be a fractal counterpart of the Gutzwiller, Chazarain and the Selberg trace formulas in this context, connecting the length spectrum of the snowflake billiard\footnote{Or, more accurately, the collection of (suitable equivalence classes of) periodic orbits of the billiard flow of $KS$ (i.e., essentially equivalently, of the geodesic flow of $\mathcal{S}_{KS}$; see part (iii) of Conjecture~\ref{conj:fractalbilliardFractalSurface}).}  and the eigenvalue \emph{(}or frequency\emph{)} spectrum of the corresponding snowflake fractal drum. Moreover, can one address the same problem for other fractal billiards \emph{(}once they have been properly defined\emph{)} and the associated fractal drums?
\end{openproblem}

We hope that the study of such open problems and conjectures, aided by a suitable combination of computer experiments and theoretical investigations, will enable us in the future to better understand the elusive nature of the Koch snowflake billiard and, eventually, of a variety of other fractal billiards.

\vspace{5 mm}
\emph{Added note.} For a version of the paper with crisper pictures, please email the second author.  Such a file will, however, be considerably larger and require between thirty and sixty minutes to download over a standard 56Kbps connection.

\vspace{5 mm}
\emph{Acknowledgements}. We wish to thank Pascal Hubert for his helpful comments on a preliminary version of this paper.

\bibliographystyle{amsalpha}

\end{document}